\documentclass[conference]{IEEEtran}
\usepackage[utf8]{inputenc}
\usepackage{amsmath}
\usepackage{amssymb}
\usepackage{pgfplots}
\usepackage{algpseudocode}
\usepackage{algorithm}
\usetikzlibrary{patterns}
\usepackage{booktabs}
\usepackage{threeparttable}

\usepackage{tikz}
\usetikzlibrary{calc}
\tikzset{generator/.style={
  circle,
  draw,
  append after command={
    \pgfextra{
    \draw
      ($(\tikzlastnode.center)!0.5!(\tikzlastnode.west)$)
       arc[start angle=180,end angle=0,radius=0.425ex] 
      (\tikzlastnode.center)
       arc[start angle=180,end angle=360,radius=0.425ex]
      ($(\tikzlastnode.center)!0.5!(\tikzlastnode.east)$) 
    ;
    }
  },
  scale=1.5,
 }
}
 \tikzset{reservoir/.style ={
  append after command={
    \pgfextra{
    \fill[rounded corners,gray!50]  ($(\tikzlastnode.north west) !0.1!(\tikzlastnode.center)-(0,1ex)$)--($(\tikzlastnode.south west) !0.1!(\tikzlastnode.center)$)--($(\tikzlastnode.south east) !0.1!(\tikzlastnode.center)$)--($(\tikzlastnode.north east) !0.1!(\tikzlastnode.center)-(0,1ex)$);
    \draw[rounded corners](\tikzlastnode.north west)--($(\tikzlastnode.north west) !0.1!(\tikzlastnode.center)$)--($(\tikzlastnode.south west) !0.1!(\tikzlastnode.center)$)--($(\tikzlastnode.south east) !0.1!(\tikzlastnode.center)$)--($(\tikzlastnode.north east) !0.1!(\tikzlastnode.center)$)--(\tikzlastnode.north east);
    \draw (\tikzlastnode.center)node{#1};
    }
  },
  scale=3.0,
 }
}

\newcommand{\NH}{\mathcal{H}} 
\newcommand{\Hindex}{h} 
\newcommand{\NT}{\mathcal{G}} 
\newcommand{\Tindex}{i} 
\newcommand{\Time}{\mathcal{T}} 
\newcommand{\Timeindex}{t} 
\newcommand{\Hydromatrix}{\mathbf{X} }
\newcommand{\Thermalmatrix}{\mathbf{Y} }

\begin{document}

\title{Inexactness of the Hydro-Thermal Coordination Semidefinite Relaxation}

\author{\IEEEauthorblockN{M.~Paredes}
\IEEEauthorblockA{IBM Research\\
São Paulo, Brazil 04007--900\\
Email: mparedes@br.ibm.com}
\and
\IEEEauthorblockN{L.~S.~A.~Martins}
\IEEEauthorblockA{IBM Research\\
São Paulo, Brazil 04007--900\\
Email: leonardo.martins@br.ibm.com}
}


\maketitle

\begin{abstract}
Hydro-thermal coordination is the problem of determining the optimal economic dispatch of hydro and thermal power plants over time.
The physics of hydroelectricity generation is commonly simplified in the literature to account for its fundamentally nonlinear nature.
Advances in convex relaxation theory have allowed the advent of Shor's semidefinite programming (SDP) relaxations of quadratic models of the problem.
This paper shows how a recently published SDP relaxation is only exact if a very strict condition regarding turbine efficiency is observed, failing otherwise.
It further proposes the use of a set of convex envelopes as a strategy to successfully obtain a stricter lower bound of the optimal solution.
This strategy is combined with a standard iterative convex-concave procedure to recover a stationary point of the original non-convex problem.
\end{abstract}


%
\IEEEpeerreviewmaketitle

\section{Introduction}
\label{sec:introduction}
Hydroelectric power generation derives from total physical work available from water elevated by dams. It is commonly expressed as an increasing function of net head ($h_n$), and turbine-released, or (equivalently) discharged water ($q$) for given turbine ($\eta_{T}$) and generator ($\eta_{G}$) efficiencies, such that
\begin{equation}
\label{eq:introduction:production}%
P_{h} =  g\cdot\rho\cdot\eta_G\cdot\eta_T\cdot q\cdot h_n,
\end{equation}
where $g$ and $\rho$ are constants representing gravity acceleration, and water density, respectively. Net head accounts for the difference between forebay ($h_{b}$) and tailwater ($h_{t}$) elevations, as well as losses due to hydraulic load ($h_l$) and atmospheric pressure differences ($h_a$)
\begin{equation}
h_{n}(\cdot) = h_{b}(v) - h_{t}(q) - h_{l}(q) - h_{a}(v).
\end{equation}
Hydraulic load losses are commonly formulated as a convex quadratic function~\cite{982207} over $q$, and losses due to atmospheric pressure are more prominent as $h_b - h_t$ increases.
Forebay elevation is a function of volume ($v$) of water in the reservoir.
Analogously, tailwater elevation is a function of water discharge.
Alternatively, it could also be a function of spillage~\cite{1602016}.
Both functions are strictly increasing on their variables if a three-dimensional geometric reservoir model is considered~\cite{VIEIRA2015781}.
If a cubic geometry is assumed then $h_b(v)$ and $h_t(q)$ are described by linear functions.
Otherwise, if trapezoidal geometries are assumed then higher order polynomials are necessary in order to represent variable head.

Nonlinear productivity with respect to water discharge in the hydroelectric power generation function, along with the intrinsic uncertainty with respect to future water availability, comprise two of the major computational challenges of the hydro-thermal coordination (HTC) problem.
Stochastic approaches not uncommonly assume a two-dimensional representation of reservoirs, i.e. constant $h_n$ and, therefore, constant productivity, resulting in linear functions of discharge~\cite{Pereira1991}.
On the other hand, deterministic models of the HTC problem resort to nonlinear formulations of the hydroelectric generation function, from second-order concave approximations~\cite{4113907,260860} to higher-order non-convex polynomial representations~\cite{Martins2014}.

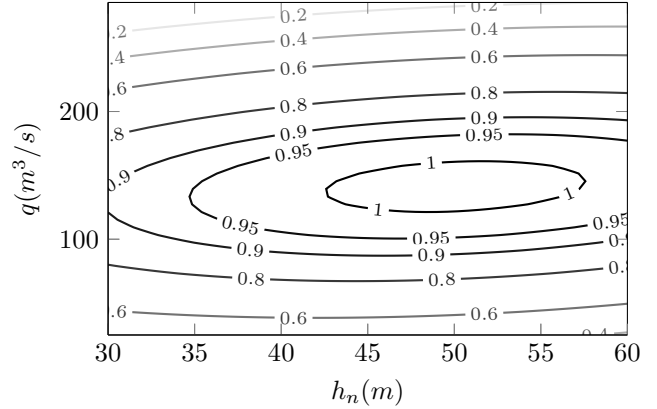
\begin{figure}
\begin{tikzpicture}
\begin{axis}[,xlabel = $h_{n}(m)$
     , ylabel = $q (m^{3}/s)$, domain=30:60,domain y =25:320 , view={0}{90},colormap={examplemap}{rgb=(0.9,0.9,0.9) rgb=(0,0,0)},width=8.5cm,height =6cm ]
 \addplot3[contour gnuplot={levels={0.2,0.4,0.6,0.8,0.9,0.95,1}},samples=50,thick] {-0.21311 + 0.022762*x+0.0093291*y+0.0000451*x*y-0.000291*x^2-0.000041*y^2   };
 \end{axis}
\end{tikzpicture}
\caption{A hypothetical hydroelectric turbine efficiency curve~\cite{466476}.}
\label{fig:efficiency}%
\end{figure}

General efficiency of a hydraulic turbine is defined as the ratio of power delivered to the shaft to the power available in the moving water.
Maximum hydraulic efficiency is specified at design time for given reference values of net head and discharge.
It is commonly described in the literature as a normalized concave quadratic function~\cite{4275241,466476} of $h_n$, and $q$,
\begin{equation}
\eta_{T}(\cdot) = e_{0} + e_{h}h_{n}+e_{q}q+e_{hq}h_{n} q+e_{hh}h_{n}^2+e_{h}q^2,
\end{equation}
as illustrated in~\figurename~\ref{fig:efficiency}.

As far as turbine and generator efficiencies are concerned, other sorts of simplifications are commonly proposed and have to do with the time resolution under consideration.
In close to real-time operation decision-making, it is advisable to take into account turbine efficiency~\cite{Finardi2006,Yeh2013}.
As the time horizon of the HTC problem increases, lower time resolutions are commonly considered, and thus the use of average efficiency is a common model assumption~\cite{VIEIRA2015781,Martins2014}.

Recently, Yunan et al.~\cite{6514678} proposed a Shor's' semidefinite relaxation of the short-term HTC problem for global solution with quadratic formulations of Eq.~\eqref{eq:introduction:production}.
In this paper we show that the relaxation suggested in~\cite{6514678} is only exact if concavity of~\eqref{eq:introduction:production} is assumed.
Moreover, we show how such hypothesis represents a strong assumption regarding turbine efficiency not reasonable in the short-run operation of a hydro plant.
Additionally, we propose the use of McCormick convex envelopes as a strategy to tighten the relaxation, in combination with a standard iterative convex-concave procedure to recover a stationary point of the original non-convex problem.

\section{Hydro-thermal coordination}%
\label{sec:htc}

Let forebay and tailwater elevations be described by linear functions of $v$ and $q$, respectively, thus implying cubic geometries of the respective reservoirs, then we have:
\begin{align}
  \label{eq:htc:hb}%
  h_{b} = & h_{b0}+h_{b1}\cdot v, \text{ and}\\
  \label{eq:htc:ht}%
  h_{t} = & h_{t0}+h_{t1}\cdot q.
\end{align}
For the sake of consistency with the related literature and the purpose of the present work, losses by hydraulic load and atmospheric pressure differences are assumed constant herein.
If generator efficiency is also assumed constant then we can write $\kappa = g\cdot\rho\cdot\eta_G$, and therefore:
\begin{displaymath}
P_{h} = \kappa\cdot\eta_T\cdot q\cdot(h_{b0}+h_{b1}\cdot v - h_{t0} - h_{t1}\cdot q - \tilde{h}_l - \tilde{h}_a).
\end{displaymath}

If, as it is commonly assumed in the longer-term HTC literature, $\eta_T$ is considered constant, e.g.\ historical mean value, then the hydroelectric power generation function can be expressed as a quadratic function
\begin{equation}
  \label{eq:prodfunc}%
  P_h = \varepsilon_q \cdot q + \varepsilon_{qq} \cdot q^{2} + \varepsilon_{qv} \cdot v \cdot q
\end{equation}
where $\varepsilon_{q}>0$,  $\varepsilon_{qq}<0$, $\varepsilon_{qv} \geqslant 0$, and $\varepsilon_{vv}=0$. Therefore, Eq.~\eqref{eq:prodfunc} is concave with respect to $q$ and indefinite with respect to $v$. In matrix form it is expressed as follows:
\begin{equation}
P_{h} ={ \mathbf{x}^{\Timeindex,\Hindex}} ^{\intercal} \widehat{\mathbf{H}}^{\Hindex}\mathbf{x}^{\Timeindex,\Hindex} +{\mathbf{e}^h}^{\intercal}\mathbf{x}^{\Timeindex,\Hindex},
\end{equation}
where
\begin{align*}
  \mathbf{x}^{\Timeindex, \Hindex} &= \left[ \begin{matrix}
  v_{\Timeindex,\Hindex}\\
  q_{\Timeindex,\Hindex}\\
  \end{matrix}\right],\\
  \widehat{\mathbf{H}}^{\Hindex} &= \left[ \begin{matrix}
 0 & \varepsilon_{qv}/2\\
\varepsilon_{qv}/2 & \varepsilon_{qq}
\end{matrix}  \right], \text{ and}\\
  \mathbf{e}^{h} &= \left[ \begin{matrix}
 0\\
 \varepsilon_{q}\\
 \end{matrix} \right].
\end{align*}

On the other hand, if $\eta_T$ is considered variable with respect to either $v$ or $q$, higher order terms arise in the hydroelectric generation function.
If, for example, we assume $\eta_T$ to be constant with respect to $q$, and, say, to vary linearly (either increasing or decreasing) with $v$, then a third-order term $vq^2$ arises, which should be ignored if one seeks to constrain~\eqref{eq:prodfunc} to the quadratic order.
Thus in this case, the concavity of $P_h$ is not necessarily defined as shown in~\cite{4113907,260860,6514678}, unless one assumes that the hydraulic turbine efficiency monotonically decreases with head, i.e. $\varepsilon_{vv} < 0$, since~\eqref{eq:htc:hb} is increasing.
This should be true only for restricted values of $v$, which is not a reasonable assumption of more general acceptability, except in cases when the hydro plant is operating with high enough head, i.e.~higher than the reference head.

The HTC problem is then formulated as that of minimizing the variable costs associated with thermoelectric power generation, subject to equations representing power balance between generation and load, and mass conservation of water, as well as inequalities representing engineering constraints, i.e.\ limits on water storage, discharge and spillage, and power output.
A Shor's semidefinite relaxation of the HTC problem is presented in~\cite{6514678} for a quadratic formulation of the hydroelectric generation function, that makes no assumption regarding the value of $\varepsilon_{vv}$.
In such relaxation we have:
\begin{align}
\label{eq:htc:X}%
\Hydromatrix^{\Timeindex, \Hindex} &= \left[ \begin{matrix}
\widehat{\Hydromatrix}^{\Timeindex, \Hindex} & \mathbf{x}^{\Timeindex, \Hindex} \\
{ \mathbf{x}^{\Timeindex, \Hindex}}^{\intercal} &  1
\end{matrix} \right], \text{ and}\\
\label{eq:htc:Y}%
\Thermalmatrix^{\Timeindex,\Tindex} &= \left[p_{\Timeindex,\Tindex} \quad 1 \right]^{\intercal}\left[p_{\Timeindex,\Tindex} \quad 1 \right],
\end{align}
where $p_{\Tindex,\Timeindex}$ represents thermoelectric power generation, such that the relaxed HTC problem is formulated as follows:
\begin{alignat}{4}
&& \min_{\Hydromatrix^{\Timeindex,\Hindex}, \Thermalmatrix^{\Timeindex,\Tindex}, \mathbf{s}} \sum_{\Tindex \in \NT,\Timeindex \in \Time} \mathbf{C}^{\Timeindex,\Tindex}\bullet\mathbf{\Thermalmatrix }^{\Timeindex,\Tindex} &  \label{eq:cost}\\
\text{s.t. } & &\sum_{\Hindex \in \NH}\mathbf{H}^{\Hindex} \bullet \Hydromatrix^{\Timeindex, \Hindex} +\sum_{\Tindex \in \NT}\mathbf{P} \bullet \Thermalmatrix^{\Timeindex,\Tindex} & \geqslant  d_{\Timeindex} & {\scriptstyle\forall  \Timeindex \in \Time} \label{eq:powerbalance}\\
&&\left( \theta_{t} \mathbf{V}+\mathbf{Q} \right)\bullet \Hydromatrix^{\Timeindex,\Hindex}-\theta_{t}\mathbf{V}\bullet \Hydromatrix^{\Timeindex-1,\Hindex} & \nonumber \\
&&+ s_{\Timeindex,\Hindex} -\sum_{\widehat{\Hindex}\in \Psi_{h}} \left( \mathbf{Q}\bullet \Hydromatrix^{\Timeindex,\widehat{\Hindex}}+ s_{\Timeindex,\widehat{\Hindex}}  \right) &=  e_{\Hindex,\Timeindex}  & {\scriptstyle\forall \Hindex \in \NH ; \Timeindex \in \Time } \label{eq:waterbalance} \\
 &&\underline{v}_{\Hindex} \leqslant \mathbf{V} \bullet \Hydromatrix^{\Timeindex,\Hindex} & \leqslant  \overline{v}_{h} & {\scriptstyle\forall \Hindex \in \NH; \Timeindex \in \Time} \label{eq:volumelimits} \\
&&\underline{q}_{\Hindex} \leqslant \mathbf{Q} \bullet \Hydromatrix^{\Timeindex,\Hindex}& \leqslant  \overline{q}_{\Hindex} & {\scriptstyle\forall \Hindex \in \NH
; \Timeindex \in \Time } \label{eq:turbinatedlimits}\\
  & &\underline{p}_{\Tindex} \leqslant \mathbf{P} \bullet \Thermalmatrix^{\Timeindex,\Tindex} & \leqslant \overline{p}_{\Tindex} & {\scriptstyle\forall \Tindex \in \NT  ; \Timeindex \in \Time }  \label{eq:powerlimits}\\
 &&\Hydromatrix^{\Timeindex, \Hindex} &\succeq 0 &  {\scriptstyle\forall \Hindex \in \NH
; \Timeindex \in \Time } \label{eq:hydrosemidefinite}\\
 &&\Thermalmatrix^{\Timeindex,\Tindex} &\succeq  0 &{\scriptstyle\forall \Tindex \in \NT  ; \Timeindex \in \Time } \label{eq:thermalsemidefinite}\\
 &&s_{\Timeindex,\Hindex} &\geqslant 0 & {\scriptstyle\forall \Hindex \in \NH ; \Timeindex \in \Time } \label{eq:snonnegative}
\end{alignat}
where $\mathbf{C}^{\Timeindex,\Tindex}$ represents the variable costs of a thermal plant $\Tindex$ at time $\Timeindex$, typically formulated as a convex function. The hydroelectric generation function is represented by $\mathbf{H}^{\Hindex}$ for a hydro plant $\Hindex$, such that:
\begin{equation}
\mathbf{H}^{\Hindex} =
\left[ \begin{matrix}
\widehat{\mathbf{H}}^{\Hindex}& \mathbf{e}^{\Hindex}\\
{\mathbf{e}^{\Hindex} }^{\intercal} & 0
\end{matrix} \right].
\end{equation}
The set of hydro plants immediately upstream of $\Hindex$ is represented by $\Psi_{h}$. A time-dependent volume-to-flow conversion coefficient is given by $\theta_{t}$. Water storage volume and discharge are represented by $\mathbf{V}$ and $\mathbf{Q}$, as well as thermoelectric power generation is equivalently represented by the Frobenius product $\mathbf{P} \bullet \Thermalmatrix^{\Timeindex,\Tindex}$, with limits defined in~\eqref{eq:powerlimits}, such that:
\begin{displaymath}
  \mathbf{V} = \left[
  \begin{matrix}
  0 & 0 & \tfrac{1}{2}\\
  0 & 0 & 0\\
  \tfrac{1}{2} & 0 & 0\\
\end{matrix} \right], \;
  \mathbf{Q} = \left[
  \begin{matrix}
  0 & 0 &0\\
  0 & 0 &  \tfrac{1}{2}\\
  0 & \tfrac{1}{2} & 0\\
  \end{matrix} \right], \text{ and }
  \mathbf{P} = \left[
  \begin{matrix}
   0 &  \tfrac{1}{2}\\
  \tfrac{1}{2} & 0\\
  \end{matrix} \right].
\end{displaymath}

In~\eqref{eq:waterbalance}, mass conservation of water is formulated as a linear algebraic system describing reservoir cascades and temporal coupling, where $s_{\Timeindex,\Hindex}$ is a variable representing spillage, and $e_{\Timeindex,\Hindex}$ is given as inflow.
Limits on volume and discharge are described in~\eqref{eq:volumelimits}, and~\eqref{eq:turbinatedlimits}, respectively.
Initial ($v_{0,\Hindex}$) and final target ($v_{T,\Hindex}$) volumes are given as boundary conditions of the problem, and respectively represented in~\eqref{eq:waterbalance}, and~\eqref{eq:volumelimits}.
As a consequence of~\eqref{eq:htc:X} and~\eqref{eq:htc:Y}, semidefiniteness constraints on $\Hydromatrix^{\Timeindex, \Hindex}$ and $\Thermalmatrix^{\Timeindex, \Tindex}$ in~\eqref{eq:hydrosemidefinite}, and~\eqref{eq:thermalsemidefinite}, respectively, whereas their respective rank-1 constraints are relaxed.

\section{Exactness of the semidefinite relaxation}%
\label{sec:relaxation}
Although power balance~\eqref{eq:powerbalance} between load, and total hydro and thermoelectric generation is a constraint that must be strictly~\cite{6514678} observed, it can be formulated as an inequality active in an optimal solution, if~\eqref{eq:cost} is a monotonically increasing function, and the following condition holds:
\begin{equation}
d_{\Timeindex} \geqslant \sum_{\Tindex \in \NT} \underline{p}_{\Tindex} + \max_{\Hydromatrix^{\Timeindex, \Hindex} \in \Omega} \sum_{\Hindex \in \NH}\mathbf{H}^{\Hindex} \bullet \Hydromatrix^{\Timeindex, \Hindex} \; \forall \Timeindex \in \Time \label{eq:maxgen}
\end{equation}
where $\Omega = \left\{ \mathbf{X}^{\Timeindex,\Hindex} : \mathbf{X}^{\Timeindex,\Hindex} \in \eqref{eq:waterbalance}, \eqref{eq:volumelimits}, \eqref{eq:turbinatedlimits} \text{ and } \eqref{eq:hydrosemidefinite} \right\}$.
In other words, condition~\eqref{eq:maxgen} establishes the reasonable assumption that, as long as load demand cannot be met exclusively with hydroelectric power, power balance equations can be exactly relaxed into inequalities in a stationary point.

Unless empirically defined by means of statistical regression with a concavity constraint, and subject to overestimation errors since higher order negative terms are dropped, it is physically reasonable to observe that, in a quadratic formulation of the hydroelectric power generation function with constant turbine efficiency (e.g. average), coefficients $\varepsilon_{0}$, $\varepsilon_{v}$, and  $\varepsilon_{vv}$ must be zero.
This results in indefiniteness of the function, as its respective eigenvalues of $\widehat{\mathbf{H}}^{\Hindex}$ are given by:
\begin{equation}
\lambda_{1},\lambda_{2} = \frac{\varepsilon_{qq} \pm \sqrt{4\varepsilon_{qv}^{2} + \varepsilon_{qq}^{2} }}{2}
\end{equation}
If, however, $\widehat{\mathbf{H}}^{\Hindex}$ is at least negative semidefinite, such that $\varepsilon_{vv} < 0$, then $P_h$ becomes concave, and the HTC problem as formulated in~\eqref{eq:cost}--\eqref{eq:snonnegative} is convex, and thus no convex relaxation is necessary.
Moreover, as shown in lemmas~1 and~2 of~\cite{2017arXiv170307870P}, Shor's semidefinite relaxation of such quadratically constrained quadratic problems (QCQP) is exact.

Despite the purportedly generality of the relaxation proposed in~\cite{6514678}, all numerical case studies presented therein fall in such QCQP convex formulation.
Its general applicability to non-concave formulations of the hydroelectric power generation function, however, fails the general conditions for relaxation exactness of~\cite{Kim2003}, since $\mathbf{V}$, $\mathbf{Q}$, and $\mathbf{P}$ are off-diagonal nonnegative.
Furthermore, the sign definiteness conditions of Sojoudi et al.~\cite{7039733} for exact relaxation cannot be confirmed since no assumptions regarding the signs of constraint coefficients are provided by Yunan et al.~\cite{6514678}.

\section{McCormick convex envelopes}%
\label{sec:mccormick}
The implicit semidefiniteness constraints (by taking the Shur's complement on $\Hydromatrix^{\Timeindex, \Hindex}$) in Shor's relaxation constitutes a lower bound of every bilinear and quadratic term in~$\widehat{\mathbf{X}}^{\Timeindex,\Hindex}$:
\begin{equation}
  \Hydromatrix^{\Timeindex, \Hindex} \succeq 0 \rightarrow \widehat{\mathbf{X}}^{\Timeindex,\Hindex} \succeq { \mathbf{x}^{\Timeindex,\Hindex}}^{\intercal}\mathbf{x}^{\Timeindex,\Hindex}.
\end{equation}
Because the physically-derived hydroelectric power generation function is indefinite, its terms with nonnegative coefficients will require an upper bound.
These upper bounds can be obtained by using McCormick convex envelopes~\cite{McCormick1976}.
This approach is based on the relaxation of bilinear terms, whose generalizations have been proposed in~\cite{Sherali1992,Anstreicher2012}.
These are called reformulation-linearization and can use the box constraints on $v_{h,t}$ and $q_{h,t}$ to construct bounds for quadratic and bilinear terms resulting from Shor's semidefinite relaxation:
\begin{equation}
(\overline{v}_{\Hindex}-v_{\Timeindex,\Hindex}) ,(v_{\Timeindex,\Hindex}-\underline{v}_{\Hindex}), (\overline{q}_{\Hindex}-q_{\Timeindex,\Hindex}) ,(q_{\Timeindex,\Hindex}-\underline{q}_{\Hindex}) \geqslant 0. \label{eq:bounds}
\end{equation}
We can multiply some of these nonnegative differences and obtain the following inequalities:
\begin{align}
\left(\overline{v}_{\Hindex} +\underline{v}_{\Hindex} \right) v_{\Timeindex,\Hindex} -\overline{v}_{\Hindex} \underline{v}_{\Hindex} & \geqslant v_{\Hindex,\Timeindex}^{2},\\
\underline{q}_{\Hindex}v_{\Timeindex,\Hindex} +\overline{v}_{\Hindex} q_{\Timeindex,\Hindex} -\overline{v}_{\Hindex}\underline{q}_{\Hindex} & \geqslant v_{\Timeindex,\Hindex}q_{\Timeindex,\Hindex}, \text{ and}\\
\underline{v}_{\Hindex}q_{\Timeindex,\Hindex} +\overline{q}_{\Hindex} v_{\Timeindex,\Hindex} -\overline{q}_{\Hindex}\underline{v}_{\Hindex} & \geqslant v_{\Timeindex,\Hindex}q_{\Timeindex,\Hindex}.
\end{align}
It follows from the superlinear monotonicity of~\eqref{eq:cost} and the additivity of~\eqref{eq:powerbalance} that hydroelectric power generation is maximized complementarily to thermal power. Therefore, given that $\varepsilon_{qq}< 0$, and $\varepsilon_{qv} \geqslant 0$, then upper bounds on bilinear terms $q_{\Timeindex,\Hindex}v_{\Timeindex,\Hindex}$ are introduced, if $\varepsilon_{qv} > 0$, by means of the following inequalities:
\begin{align}
\left[ \begin{matrix}
 0 & -\tfrac{1}{2} & \tfrac{1}{2} \underline{q}_{\Hindex}  \\
  - \tfrac{1}{2} &  0 & \tfrac{1}{2} \overline{v}_{\Hindex} \\
   \tfrac{1}{2} \underline{q}_{\Hindex}  &\tfrac{1}{2} \overline{v}_{\Hindex}  & 0
\end{matrix} \right]  \bullet \Hydromatrix^{\Timeindex,\Hindex} &\geqslant \overline{v}_{\Hindex}\underline{q}_{\Hindex}, \text{ and}\label{eq:q1m}\\
\left[ \begin{matrix}
 0 & -\tfrac{1}{2} & \tfrac{1}{2} \overline{q}_{\Hindex}  \\
  - \tfrac{1}{2} &  0 & \tfrac{1}{2} \underline{v}_{\Hindex} \\
   \tfrac{1}{2} \overline{q}_{\Hindex}  &\tfrac{1}{2} \underline{v}_{\Hindex}  & 0
\end{matrix} \right]  \bullet \Hydromatrix^{\Timeindex,\Hindex} &\geqslant \underline{v}_{\Hindex}\overline{q}_{\Hindex}. \label{eq:q2m}
\end{align}
Analogously, $v^{2}$ is also upper bounded:
\begin{equation}\label{eq:upperv2}
\left[ \begin{matrix}
 -1 & 0 & \tfrac{1}{2} \left(\overline{v}_{\Hindex} +\underline{v}_{\Hindex} \right)   \\
  0 &  0 &  0 \\
   \tfrac{1}{2}  \left(\overline{v}_{\Hindex} +\underline{v}_{\Hindex} \right)   &0  & 0
\end{matrix} \right]  \bullet \Hydromatrix^{\Timeindex,\Hindex} \geqslant \overline{v}_{\Hindex} \underline{v}_{\Hindex}.
\end{equation}





\section{Convex-concave procedure}%
\label{sec:covexconcave}
Any mathematical program with indefinite quadratic functions in the constraint set or the objective can be expressed as a difference of convex (DC) programming problems~\cite{Lipp2016} of the following form:
\begin{align}
\min_{x} &  &f_{0}(x) - g_{0}(x)& \\
\text{s.t. } &&f_{i}(x) - g_{i}(x) &\leqslant 0 \quad i=1,\ldots, m
\end{align}
where $f_{i} : \mathrm{R}^{n} \rightarrow \mathrm{R}$ and $g_{i} : \mathrm{R}^{n} \rightarrow \mathrm{R}$ for $i=0,\ldots,m$ are convex functions.
The convex-concave procedure (CCP) is a heuristic based on DC problems to obtain a stationary point of the original non-convex problem as explained in~\cite{NIPS2009_3646}.
In this case the DC hydroelectric power generation functions are in the power balance constraints~\eqref{eq:powerbalance}, and can be expressed as a difference of differentiable quadratic functions by decomposing $-\widehat{\mathbf{H}}^{\Hindex}$ into a difference of positive definite matrices:
\begin{equation}
\widehat{\mathbf{H}}^{\Hindex} = -\widehat{\mathbf{H}}^{\Hindex}_{+} + \widehat{\mathbf{H}}_{-}^{\Hindex}.
\end{equation}

\begin{algorithm}
\begin{algorithmic}[1]
\Require Solution $\mathbf{x}^{\Timeindex, \Hindex}_{(0)}$ to problem \eqref{eq:cost}--\eqref{eq:snonnegative}, \eqref{eq:q1m}--\eqref{eq:upperv2}
\State Decompose $\widehat{\mathbf{H}}^{\Hindex} = -\widehat{\mathbf{H}}^{\Hindex}_{+} + \widehat{\mathbf{H}}_{-}^{\Hindex} $
\State Set $k \leftarrow 0$
\Repeat
\State Set $k \leftarrow k + 1$
\State Construct $\widetilde{\mathbf{H}}^{\Hindex}_{(k)}$
\State Obtain $\mathbf{x}^{\Timeindex, \Hindex}_{(k)}$ from solving \eqref{eq:cost}, \eqref{eq:QC2}, \eqref{eq:waterbalance}--\eqref{eq:snonnegative}
\Until{$ \sum_{\Tindex \in \NT,\Timeindex \in \Time}\left\| \mathbf{C}^{\Timeindex,\Tindex} \bullet \left( \mathbf{\Thermalmatrix }^{\Timeindex,\Tindex}_{(k)}-\mathbf{\Thermalmatrix }^{\Timeindex,\Tindex}_{(k-1)} \right) \right\|_{2} \leq \epsilon $}
\end{algorithmic}%
\caption{Convex-concave procedure.}
\label{alg:ccp}
\end{algorithm}

Under CCP the concave part of the indefinite quadratic function becomes affine by means of a first-order Taylor series approximation around $\mathbf{x}^{\Timeindex, \Hindex}_{(k-1)}$:
\begin{equation*}
\widetilde{\mathbf{H}}^{\Hindex}_{(k)} = \left[
\begin{matrix}
-\widehat{\mathbf{H}}^{\Hindex}_{+} & \tfrac{1}{2}  \left(  \mathbf{e}^{\Hindex} + \widehat{\mathbf{H}}_{-}^{\Hindex} {\mathbf{x}^{\Timeindex, \Hindex}_{(k-1)}} \right)\\
\tfrac{1}{2}  \left( \mathbf{e}^{\Hindex} + \widehat{\mathbf{H}}_{-}^{\Hindex} {\mathbf{x}^{\Timeindex, \Hindex}_{(k-1)}} \right)^{\intercal} & {\mathbf{x}^{\Timeindex, \Hindex}_{(k-1)}}^{\intercal}  \widehat{\mathbf{H}}_{-}^{\Hindex}\mathbf{x}^{\Timeindex, \Hindex}_{(k-1)}
\end{matrix}
\right],
\end{equation*}
where $k$ is the $k$-th iteration of Algorithm~\ref{alg:ccp}, such that the power balance constraints are iteratively reformulated:
\begin{equation}
\sum_{\Hindex \in \NH}\widetilde{\mathbf{H}}^{\Hindex}_{(k)}\bullet \Hydromatrix^{\Timeindex, \Hindex} +\sum_{\Tindex \in \NT}\mathbf{P} \bullet \Thermalmatrix^{\Timeindex,\Tindex}  \geq d_{\Timeindex}. \label{eq:QC2}
\end{equation}
Each solution $\mathbf{x}^{\Timeindex, \Hindex}_{(k)}$ is recovered as follows:
\begin{displaymath}
  x^{\Timeindex, \Hindex}_{(k),j} = \sqrt{\Hydromatrix^{\Timeindex, \Hindex}_{(k),jj}}, \quad j = 1,2.
\end{displaymath}

\section{Numerical experiments}%
\label{sec:results}
A numerical example is shown to illustrate the inexactness of the Shor's SDP relaxation of Yunan et al.~\cite{6514678} for the simple case in which an average hydraulic efficiency is considered, thus resulting in a non-concave hydroelectric production function as the one formulated in~\eqref{eq:prodfunc}.
Complete case study data are presented in Appendix~\ref{appendix1}.
The case study uses data from 5 hydro plants in the Brazilian Paranaíba river basin.
A fictitious thermal plant complements the hypothetical case study power system.
Hydro plants GH2, GH4, and GH5 are run-off-river, meaning that their reservoir volumes remain constant for all $\Timeindex$ with monthly discretization in a year span.
Boundary conditions at maximum storage volume were equally defined for each of the hydro plants.
\figurename~\ref{fig:system} depicts the transmission-unconstrained system configuration with constant $1551.4$~MW load.
Optimization was carried over in Python~3 and CVXPY~\cite{cvxpy} interfaced with the SDPA~\cite{SDPA} solver for semidefinite programming.

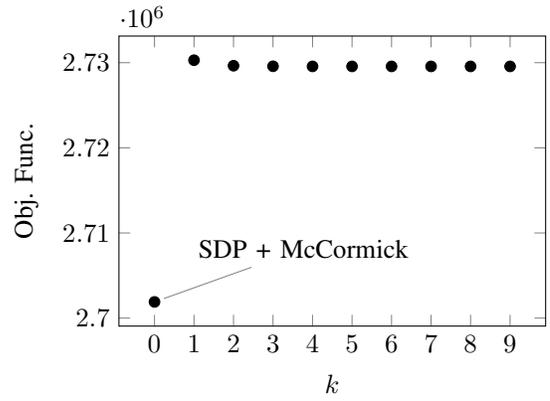
\begin{figure}
  \centering
\begin{tikzpicture}
    \begin{axis}[
        height=0.3\textwidth,
        width=0.4\textwidth,
        xlabel=$k$,
        xtick={0,...,9},
        ylabel=Obj. Func.,
]
        \addplot[only marks,mark=*,black] coordinates {(0,2701900.11389) (1,2730291.78877) (2,2729645.90622)(3,2729579.66422)(4,2729568.20284) (5,2729566.03575) (6,2729565.61775)(7,2729565.53695)(8,2729565.52129) (9,2729565.51821) }
        node[pos=0.0, pin=45:SDP + McCormick]{} ;
    \end{axis}
\end{tikzpicture}
\caption{Progress of the objective function value at each iteration.}%
\label{fig:iterations}
\end{figure}

A comparison of the lower bounds provided by the different approaches is listed in Table~\ref{tab:comparison}, along with the stationary point found by CCP.
The introduction of McCormick convex envelopes has allowed for a two orders of magnitude improvement on the objective function lower bound.
\figurename~\ref{fig:iterations} illustrates the objective function values at the end of each of the 9 CCP iterations required for convergence in about 7 seconds.
Table~\ref{tab:powerresult} lists the power generation results for each of the plants.

\begin{table}
  \centering
  \scriptsize
  \begin{threeparttable}
    \caption{Comparison between relaxation strategies.}%
    \label{tab:comparison}
    \begin{tabular}{lr}
      \toprule
      Relaxation & Obj.~Func.\\
      \midrule
      SDP\tnote{1} & 68,083.08\\
      SDP + McCormick & 2,701,900.11\\
      \midrule
      Convex-concave procedure & 2,729,565.52 \\
      \bottomrule
    \end{tabular}
    \begin{tablenotes}
      \item[1] Yunan et al.~\cite{6514678}.
    \end{tablenotes}
  \end{threeparttable}
\end{table}

\begin{table}[!t]
\caption{Power generation results.} \label{tab:powerresult}
\centering
\scriptsize
\begin{tabular}{lrrrrrr}
\toprule
Month &  GH1 & GH2 & GH3 &  GH4 &  GH5 & GT1 \\
\midrule
1&          117.49 &           98.73 &     217.53 &             304.40 &         37.06   &   776.19\\
2&           91.91 &           77.27 &     162.84 &             237.01 &         29.48 & 952.90\\
3 &           76.47 &           64.30 &     128.83 &             191.17 &         24.68 & 1,065.94\\
4 &           63.61 &           53.72 &     130.50 &             167.80 &         21.74 & 1,114.03\\
5 &           58.69 &           49.47 &     159.59 &             176.83 &         22.55 & 1,084.27\\
6 &           71.11 &           59.93 &     162.42 &             207.44 &         28.34 &   1,022.16\\
7 &          103.14 &           87.44 &     166.17 &             283.50 &         36.44 &   874.72\\
8 &          164.88 &          139.40 &     209.08 &             437.88 &         50.24 &   549.92\\
9 &          232.45 &          194.34 &     250.80 &             560.90 &         61.37 &   251.55 \\
10 &          239.82 &          200.71 &     302.70 &             577.69 &         66.59 & 163.88\\
11 &          224.86 &          188.42 &     375.03 &             580.51 &         66.02 &   116.56 \\
12 &          169.54 &          142.65 &     335.53 &             454.64 &         52.01&   397.03 \\
\bottomrule
\end{tabular}
\end{table}

\section{Conclusion}
In a recent paper by Yunan et al.~\cite{6514678} a Shor's semidefinite relaxation with global optimality was proposed for the quadratically constrained quadratic hydro-thermal coordination problem.
In this paper, however, we present an empirical analysis showing that such relaxation exactness is only known to be possible if concavity of the hydroelectric power generation function is assumed, therefore resulting, under reasonable assumptions regarding load demand, in an already fully convex problem formulation whose SDP relaxation is redundant.
Furthermore, the concavity assumption hypothesis represents a strong assumption regarding turbine efficiency only reasonable in the short-run operation of hydro plants in very limited situations.
In a numerical case study with average turbine efficiencies, and thus indefinite production functions, the use of McCormick convex envelopes was shown to provide tighter lower bounds on the objective function by orders of magnitude.
Additionally, we provide a reformulation-linearization for stationary point recovery by means of an iterative convex-concave procedure that further suggests the effectiveness of the convex envelope use.

\appendix[Case study data]%
\label{appendix1}
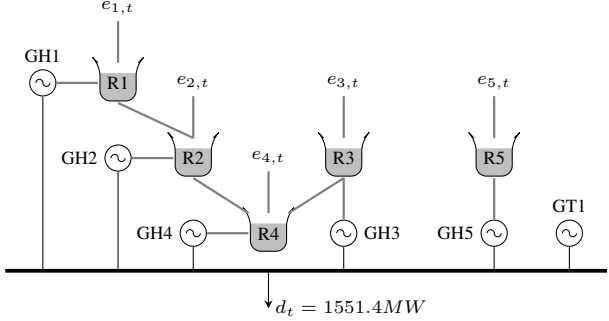
\begin{figure}[!h]%
\centering
\scriptsize
\begin{tikzpicture}
\draw
(-1,1) node [reservoir={ R1 }] (R1) {}
(-1,2) node[](a1){$e_{1,\Timeindex}$}
(0,0) node [reservoir={R2}] (R2) {}
(0,1) node[](a2){$e_{2,\Timeindex}$}
(2,0) node [reservoir={  R3}] (R3) {}
(2,1) node[](a3){$e_{3,\Timeindex}$}
(1,-1) node [reservoir={  R4}] (R4) {}
(1,0) node[](a4){$e_{4,\Timeindex}$}
(4,0) node [reservoir={  R5}] (R5) {}
(4,1) node[](a5){$e_{5,\Timeindex}$}
(-2,1) node [generator,label=above:{GH1}] (GH1){}
(-1,0) node [generator,label=left:{GH2}] (GH2){}
(2,-1) node [generator,label=right:{GH3}] (GH3){}
(0,-1) node [generator,label=left:{GH4}] (GH4){}
(4,-1) node [generator,label=left:{GH5}] (GH5){}
(5,-1) node [generator,label={ [align=center] above: GT1}] (GT1){}
;
\draw [line width=0.3mm, gray ]
 (R1.south)--(R2.north)  
(R2.south)--(R4.north west)
(R3.south)--(R4.north east)
(R1.west)--(GH1.east)
(R2.west)--(GH2.east)
(R3.south)--(GH3.north)
(R4.west)--(GH4.east)
(R5.south)--(GH5.north)
(a1.south)--(R1.north)
(a2.south)--(R2.north)
(a3.south)--(R3.north)
(a4.south)--(R4.north)
(a5.south)--(R5.north)
;
\draw
(GH1.south) --(-2,-1.5)
(GH2.south) --(-1,-1.5)
(GH3.south) --(2,-1.5)
(GH4.south) --(0,-1.5)
(GH5.south) --(4,-1.5)
(GT1.south) --(5,-1.5)
;
\draw [line width=0.5mm] (-2.5,-1.5) coordinate(bar-left) -- (5.5,-1.5)coordinate(bar-right) ;
\draw[-stealth](1,-1.5) --++(0,-0.5cm) node[right] {$d_\Timeindex=1551.4MW$}; 
\end{tikzpicture}
\caption{Case study system configuration.}%
\label{fig:system}
\end{figure}

\begin{table}[!h]
\caption{Hydro plant storage ($hm^{3}$) and discharge ($m^{3}/s$) data.}
\scriptsize
\begin{tabular}{lrrrrrrr}
\toprule
 & $\overline{v}$ & $\underline{v} $ & $\overline{q}$ & $\underline{q}$ &  $\varepsilon_{q}$ & $\varepsilon_{qq} $ & $\varepsilon_{qv}$ \\
\midrule
GH1 & 241.1 & 228.3& 483.0 & 65 & 0.297 & -3.06E-5 & 3.84E-4 \\
GH2 & 879.0 & 879.0 & 503.5 & 68 & 0.178  & -2.25E-5 & 1.50E-4 \\
GH3 & 470.0 &  1500.0 & 544.6 & 74 &  0.323 & -6.74E-5 & 1.50E-4 \\
GH4 & 460.0 & 460.0 & 2434.6 & 273 & 0.229 & -1.00E-5 & 0.00E-0 \\
GH5 & 95.3 & 95.3 & 277.2 & 62 & 0.198 & -4.08E-5 & 0.00E-0 \\
\bottomrule
\end{tabular}
\end{table}

\begin{table}[!h]
\caption{Inflows ($m^{3}/s$) and number of days/month.}%
\label{tab:hydrodata}
\scriptsize
\center
\begin{tabular}{rrrrrrr}
\toprule
 $\Timeindex$ & $e_1$ & $e_2$ & $e_3$ & $e_4$ & $e_5$ & Days \\
\midrule
1 & 228.17 &  238.51 &  300.19  & 1036.32 &  143.48 & 31\\
2 & 177.79 &  185.95 &  223.64 &  798.18  & 113.40 & 30\\
 3 & 147.57  & 154.40 &  176.42 &  639.20 &   94.58 & 31\\
 4 & 122.50  & 128.75  & 139.48 &  513.47  &  83.10 & 31\\
 5 & 112.94  & 118.47  & 120.53   & 458.40  &  86.25 & 30\\
 6 &  135.44 &  141.11 &  149.51  &  548.53  & 108.90 & 31\\
 7 &  200.30 &  208.61 &  252.84 &  889.94  & 141.00 & 30\\
 8 &  327.02  & 339.28  & 408.98 & 1524.89  & 196.73 & 31\\
 9  & 464.75  &  482.69 &   563.12 & 2134.43 &  242.70 & 31\\
 10 &  475.13  &  494.16 &  624.25 & 2231.39 &  264.60 & 28\\
 11 & 444.38 &  462.81 &  626.15 & 2176.34 &  262.20& 31\\
12 & 332.01  & 347.40 &  467.93 & 1587.95  & 204.00 & 30\\
\bottomrule
\end{tabular}
\end{table}

\begin{table}[!h]
\caption{Thermal plant data.}%
\label{tab:thermaldata}
\scriptsize
\center
\begin{tabular}{ccccc}
\toprule
$\overline{p}$ & $\underline{p}$ & $c_{0}$ & $c_{1}$ & $c_{2}$ \\
\midrule
1551.4 & 0.0 & 0.0 & 0.0 & 0.5\\
\bottomrule
\end{tabular}
\end{table}

\bibliographystyle{ieeepes}
\bibliography{references}

\end{document}